\documentclass[leqno,12pt]{article}

\usepackage{latexsym}
 \usepackage{graphicx}
\usepackage{amsmath}
\usepackage{amssymb}
\usepackage{cite}
\usepackage{color}

\usepackage[margin=1.3in, top=1.3in, bottom=1.3in]{geometry}

\newcommand{\nc}{\newcommand}
\nc{\nt}{\newtheorem}
\nt{thm}{Theorem}[section]
\nt{cor}[thm]{Corollary}
\nt{prop}[thm]{Proposition}
\nt{lem}[thm]{Lemma}
\nt{defn}[thm]{Definition}
\nt{rem}[thm]{Remark}
\nt{exa}[thm]{Example}
\nt{ass}[thm]{Assumption}
\nc{\ip}[2]{\mbox{$\langle #1,#2 \rangle$}}
\nc{\pf}{\noindent{\bf Proof\ \ }}
\nc{\finpf}{\hfill{$\Box$}\linespace}
\nc{\linespace}{\vspace{\baselineskip} \noindent}
\nc{\R}{{\bf R}}
\nc{\Rn}{{\bf R}^n}
\nc{\bx}{\bar{x}}
\nc{\e}{\epsilon}
\nc{\cl}{\mbox{\rm cl}\,}
\nc{\conv}{\mbox{\rm conv}\,}

\newenvironment{myequation}{\setcounter{equation}{\value{thm}}
   \begin{equation}}{\addtocounter{thm}{1}\end{equation}}
\newenvironment{myeqnarray}{\setcounter{equation}{\value{thm}}
    \begin{eqnarray}}{\setcounter{thm}{\value{equation}}\end{eqnarray}}

\nc{\bmye}{\begin{myequation}}
\nc{\emye}{\end{myequation}}

\begin{document}
\title{
BFGS convergence to nonsmooth minimizers of convex functions
}
\author{
 J. Guo
\thanks{ORIE, Cornell University, Ithaca, NY 14853, U.S.A.
\texttt{jg826@cornell.edu}.}
\and 
A.S. Lewis
\thanks{ORIE, Cornell University, Ithaca, NY 14853, U.S.A.
\texttt{people.orie.cornell.edu/aslewis}.
Research supported in part by National Science Foundation Grant DMS-1613996.}
}
\maketitle

\begin{abstract}
The popular BFGS quasi-Newton minimization algorithm under reasonable conditions converges globally on smooth convex functions.  This result was proved by Powell in 1976:  we consider its implications for functions that are {\em not} smooth.  In particular, an analogous convergence result holds for functions, like the Euclidean norm, that are nonsmooth at the minimizer.\end{abstract}
\medskip

\medskip

\medskip

\section{Introduction} 
The BFGS (Broyden-Fletcher-Goldfarb-Shanno) method for minimizing a smooth function has been popular for decades \cite{nocedal_wright}.  Surprisingly, however, it can also be an effective general-purpose tool for nonsmooth optimization \cite{BFGS}.  For twice continuously differentiable convex functions with compact level sets, Powell \cite{powell_global} proved global convergence of the algorithm in 1976.  By contrast, in the nonsmooth case, despite substantial computational experience, the method is supported by little theory.  Beyond one dimension, with the exception of some contrived model examples \cite{lewis-zhang}, the only previous convergence proof for the standard BFGS algorithm applied to a nonsmooth function seems to be the analysis of the two-dimensional Euclidean norm in \cite{BFGS}.

As a simple illustration, consider the nonsmooth convex function $f \colon \R^2 \to \R$ defined by 
$f(u,v) = u^2 + |v|$.  A routine implementation of the BFGS method, using a random initial point and a standard backtracking line search, invariably converges to the unique optimizer at zero.  Not surprisingly, the method of steepest descent, using the same line search, often converges to a nonoptimal point $(u,0)$ with $u \ne 0$.

\begin{figure}  \label{convergence}
\parbox{7cm}
{
\begin{center}
\includegraphics[width=7cm]{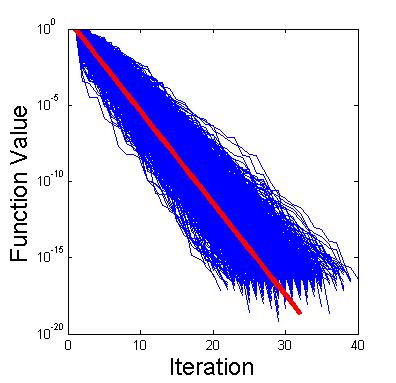}
\end{center}
}
\hfill
\parbox{7cm}
{
\begin{center}
\includegraphics[width=7cm]{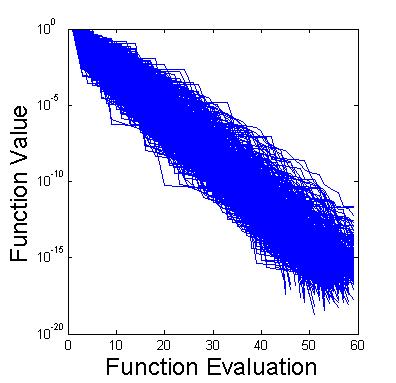}
\end{center}
}
\caption{BFGS method for $f(u,v)=u^2 + |v|$.  A thousand random starts, using inexact line search, and initial approximate Hessian $I$.  Semilog plots of function value $f(u_k,v_k)$, initially normalized.  Panel 1:  against iteration count $k$.  (Bold line plots $2^{-2k}$.)  Panel 2:  against function evaluation count, including line search.}
\end{figure}

For example, Figure \ref{convergence} plots function values for a thousand runs of BFGS against both iteration count and a count of the number of function-gradient evaluations, including those incurred in each line search.  (Precisely, the initial Hessian approximation is the identity, the weak Wolfe line search uses Armijo parameter $10^{-4}$ and Wolfe parameter $0.9$, and the initial function value is normalized to one.)  The results compellingly support convergence, and indeed suggest a linear rate:  the bold line overlaid on the first panel corresponds to the BFGS iterates $(2^{-k}, \frac{2}{5}(-1)^k 2^{-2k})$ generated by an exact line search \cite{lewis-zhang}.  However, even for this very simple example, a general convergence result does not seem easy.

Nonetheless, Powell's theory does have consequences even in the nonsmooth case.  Loosely speaking, we prove, at least under a strict-convexity-like assumption, that global convergence can only fail for the BFGS method if a subsequence of the iterates converges to a nonsmooth point.  For example, for the function 
$f(u,v)=u^2 + |v|$, BFGS iterates cannot remain a uniform distance away from the line $v=0$.  While intuitive --- a successful smooth algorithm should somehow detect nonsmoothness --- this result is also reassuring, and in fact suffices to prove convergence on some interesting examples.  An analogous technique proves convergence for the Euclidean norm on $\Rn$, generalizing the result for $n=2$ in \cite{BFGS}.

\section{BFGS sequences}
Given a set $U \subset \Rn$, we consider the BFGS method for minimizing a possibly nonsmooth function 
$f \colon U \to \R$.  We call a sequence $(x_k)$ in $U$ ``BFGS'' if the BFGS method could generate it using a line search satisfying the Armijo and weak Wolfe conditions.  More precisely, we make the following definition.

\begin{defn}
{\rm
A sequence
$(x_k)$ is a {\em BFGS sequence} for the function $f$ if $f$ is differentiable at each iterate $x_k$ with nonzero gradient $\nabla f(x_k)$, and there exist parameters $\mu <  \nu$ in the interval $(0,1)$ and an $n$-by-$n$ positive definite matrix $H_0$  such that the vectors
\[
s_k = x_{k+1} - x_k  ~~\mbox{and}~~ y_k = \nabla f(x_{k+1}) - \nabla f(x_k)
\]
and the matrices defined recursively by
\bmye \label{bfgs}
V_k = I - \frac{s_k y_k^T}{s_k^T y_k}  ~~\mbox{and}~~ 
H_{k+1} = V_k H_k V_k^T + \frac{s_k s_k^T}{s_k^T y_k}
\emye
satisfy 
\begin{myeqnarray}
H_k \nabla f(x_k) & \in & -\R_+ s_k \label{newton} \\
f(x_{k+1}) & \le & f(x_k) + \mu \nabla f(x_k)^T s_k \label{armijo} \\
\nabla f(x_{k+1})^T s_k & \ge & \nu \nabla f(x_k)^T s_k \label{wolfe}
\end{myeqnarray}
for $k=0,1,2,\ldots$.
}
\end{defn}

Notice that this property is independent of any particular line search algorithm used to generate the sequence $(x_k)$:  it depends only on the sequences of functions values $f(x_k)$ and gradients 
$\nabla f(x_k)$.  Conceptually, in the definition, the matrices $H_k$ are approximate inverse Hessians for the function $f$ at the iterate $x_k$\/:  the equations (\ref{bfgs}) define the BFGS quasi-Newton update and the inclusion (\ref{newton}) expresses the fact that the step $s_k$ is in the corresponding approximate Newton direction.  The inequalities (\ref{armijo}) and (\ref{wolfe}) are the Armijo and weak Wolf line search conditions respectively, with parameters 
$\mu$ and $\nu$ respectively.  By a simple and standard induction argument, they imply that the property $s_k^T y_k > 0$ then holds for all 
$k$, ensuring the matrices $H_k$ are well-defined and positive definite, and that the function values $f(x_k)$ decrease strictly.  An implementation of the BFGS method for a convex function $f$ using a standard backtracking line search will generate a BFGS sequence of iterates, assuming that those iterates stay in the set $U$ and that the method never encounters a nonsmooth or critical point.

\subsubsection*{Example:  a simple nonsmooth function}
Consider the function $f \colon \R^2 \to \R$ defined by $f(u,v) = u^2 + |v|$.  (We abuse notation slightly and identify the vector $[u~v]^T \in \R^2$ with the point $(u,v)$.)  Then the sequence 
in $\R^2$ defined by
\[
\Big( 2^{-k} ~,~\frac{2}{5}(-1)^k 2^{-2k} \Big) ~~~ (k=0,1,2,\ldots)
\]
is a BFGS sequence, as observed in \cite[Prop 3.2]{lewis-zhang}.  Specifically, if we define a matrix
\[
H_0 = 
\left[
\begin{array}{cc}
\frac{1}{4} & 0 \\
0 & \frac{1}{2}
\end{array}
\right],
\]
then the the definition of a BFGS sequence holds for any parameter values \mbox{$\mu \in \!(0,0.7]$} and 
$\nu \in (\mu,1)$.  In this example, the ``exact'' line search property 
\mbox{$\nabla f(x_{k+1})^T s_k = 0$} holds for all $k$, and the approximate inverse Hessians are
\[
H_1 = 
\left[
\begin{array}{cc}
\frac{1}{2} & 0 \\
0 & \frac{1}{4}
\end{array}
\right], ~~~
H_k = 
\frac{1}{6}
\left[
\begin{array}{cc}
5 & (-1)^k 2^{1-k} \\
(-1)^k 2^{1-k} & 2^{3-2k}
\end{array}
\right]
~~~ (k>1).
\]

\subsubsection*{Example:  the Euclidean norm}
Consider the function $f = \|\cdot\|$ on $\R^2$.  Beginning with the initial vector $[1~0]^T$, generate a sequence of vectors by, at each iteration, rotating clockwise through an angle of $\frac{\pi}{3}$ and shrinking by a factor $\frac{1}{2}$.  The result is a BFGS sequence for $f$, as observed in \cite{BFGS}.  Specifically, if we define a matrix
\[
H_0 = 
\left[
\begin{array}{cc}
3 & -\sqrt{3} \\
-\sqrt{3} & 3
\end{array}
\right],
\]
then the the definition of a BFGS sequence holds for any parameter values
$\mu \in \big(0,\frac{2}{3}\big]$ and any $\nu \in (\mu,1)$.  Again, the exact line search property 
$\nabla f(x_{k+1})^T s_k = 0$ holds for all $k$.   In this case the approximate inverse Hessians have eigenvalues behaving asymptotically like $2^{-k}(3 \pm \sqrt{3})$ (see \cite{BFGS}).

\section{Main result}
The following theorem captures a key global convergence property of the BFGS method.
 
\begin{thm}[Powell, 1976] \label{thm1}
Consider an open convex set $U \subset \Rn$ containing a BFGS sequence $(x_k)$ for a convex function 
$f \colon U \to \R$.  Assume that the level set $\{x \in U : f (x) \le f(x_0)\}$ is bounded, and that
\bmye \label{continuous}
\mbox{$\nabla^2 f$ is continuous throughout $U$.}
\emye
Then the sequence of function values $f(x_k)$ converges to $\min f$.
\end{thm}

Among the assumptions in Powell's theorem, at least for dimension $n>2$ (see \cite{powell2var}), convexity is central.  Although the BFGS method works well in practice on general smooth functions \cite{nocedal_wright}, nonconvex counterexamples are known where convergence fails:  in particular, \cite{Dai2013} presents a bounded but nonconvergent BFGS sequence for a polynomial $f \colon \R^4 \to \R$.  In the general {\em convex} case, on the other hand, whether the smoothness assumption (\ref{continuous}) can be weakened seems unclear.

We present here a result analogous to Powell's theorem.  We modify the assumptions, strengthening the convexity assumption but weakening the smoothness requirement (\ref{continuous}).  Similar results to the one below hold for many common minimization algorithms possessing suitable global convergence properties in the smooth case.  Such algorithms generate sequences of iterates $x_k$ characterized by certain properties of the function values $f(x_k)$ and gradients $\nabla f(x_k)$ (for $k=0,1,2,\ldots$), analogous to the definition of a BFGS sequence.  Providing the algorithm generates function values $f(x_k)$ that must decrease to the minimum value $\min f$ for any convex function whose level sets are bounded and whose Hessian is continuous and positive definite throughout those level sets, exactly the same proof technique applies.  Examples of such algorithms include standard versions of steepest descent \cite{nocedal_wright}, coordinate descent (see for example \cite{luo-tseng-coordinate}), and conjugate gradient methods (see for example \cite{gilbert-nocedal}).  Here we concentrate on BFGS because, in striking contrast to these methods, the BFGS method works well in practice on nonsmooth functions \cite{BFGS}.

\begin{thm} \label{thm2}
Powell's Theorem also holds with the smoothness assumption (\ref{continuous}) replaced by the following assumption:
\bmye \label{modified}
\left\{
\begin{array}{l}
\mbox{$\nabla^2 f$ is positive-definite and continuous throughout}  \\
\mbox{an open set $V \subset U$ containing the set $\cl\!(x_k)$ and}  \\
\mbox{satisfying $\inf_V f = \min f$.}
\end{array}
\right.
\emye
\end{thm}

\pf
We consider an open convex set $U \subset \Rn$ containing a 
BFGS sequence $(x_k)$ for a convex function $f \colon U \to \R$ satisfying assumption (\ref{modified}).
We further assume that the level set $\{x \in U : f (x) \le f(x_0)\}$ is bounded, and our aim is to prove that the sequence of function values $f(x_k)$ converges to $\min f$. 

Assume first that the theorem is true in the special case when $U = \Rn$ and the complement $V^c$ is bounded.  We then deduce the general case as follows.  Note by assumption, that the function $f$ is not constant, so by convexity there exists a point $\bx \in U$ with $f(\bx) > f(x_0)$.  Convexity also ensures that $f$ is 
$L$-Lipschitz on the nonempty compact convex set 
\[
K ~=~ \{ x \in U : f(x) \le f(\bx) \}, 
\]
for some constant $L > 0$.  Hence there exists  a convex Lipschitz function $\hat f \colon \Rn \to \R$ agreeing with $f$ on $K$, specifically the Lipschitz regularization defined by
\[
\hat f(y) ~=~ \min_{x \in K}\{f(x) + L\|y-x\|\} ~~~ (y \in \Rn).
\]
Now, for any sufficiently large $\beta \in \R$, the convex function $\tilde f \colon \Rn \to \R$ defined by 
\[
\tilde f (x) = \max \Big\{ \hat f(x), \frac{1}{2}\|x\|^2 - \beta \Big\}   ~~~ (x \in \Rn)
\]
also agrees with $f$ on $K$.  The Hessian of $\tilde f$ is just the identity throughout the open set  
\[
W ~=~ \Big\{x : \hat f(x) < \frac{1}{2}\|x\|^2 - \beta \Big\}.
\]
Furthermore, this set has bounded complement, and therefore so does the open set
\[
\tilde V ~=~ W ~\cup~ \{ x \in V : f(x) < f(\bx) \}.
\]
Now notice that $(x_k)$ is also a BFGS sequence for the function $\tilde f$, and all the assumptions of the theorem hold with $f$ replaced by $\tilde f$, $U$ replaced by $\Rn$, and $V$ replaced $\tilde V$.  Applying the special case of the theorem, we deduce 
\[
f(x_k) = \tilde f(x_k) \to \min\tilde f = \min f, 
\]
as required.

We can therefore concentrate on the special case when $U = \Rn$ and the set $N = V^c$ is compact.  We can assume $N$ is nonempty, since otherwise the result follows immediately from Powell's Theorem.  The convex function $f$  is then continuous throughout $\Rn$.  It is not constant, and hence is unbounded above.  Furthermore, by assumption, the initial point $x_0$ is not a minimizer, so all the level sets 
\mbox{$\{ x : f(x) \le \alpha \}$} are compact.  Since $N$ is compact and $f$ is continuous, we can fix a constant $\alpha > f(x_0)$ satisfying $\alpha > \max_N f$.

Since the values $f(x_k)$ are decreasing, the sequence $(x_k)$ is bounded and hence the closure 
$\cl\!(x_k)$ is compact.   For all sufficiently small $\e > 0$, we then have
\[
\cl\!(x_k) \cap (N+2\e B) = \emptyset ~~\mbox{and}~~ \max_{N+2\e B} f < \alpha,
\]
where $B$ denotes the closed unit ball in $\Rn$.  The distance function $d_N \colon \Rn \to \R$ defined by
$d_N(x) = \min_N \|\cdot - x\|$ (for $x \in \Rn$) is continuous, so the set
\[
\Omega_{\e} ~=~ \{x : d_N(x) \ge 2\e ~\mbox{and}~ f(x) \le \alpha \}
\]
is compact, and is contained in the open set $\{ x : d_N(x) > \e \}$.  On this open set, the function $f$ is convex, in the sense of \cite{yan}, and $C^{(2)}$ with positive-definite Hessian.  Hence, by \cite[Theorem 3.2]{yan}, there exists a $C^{(2)}$ convex function $f_{\e}$ on a convex open neighborhood 
$U_{\e}$ of the convex hull $\conv \Omega_{\e}$ agreeing with $f$ on $\Omega_{\e}$.  Our choice of $\e$ ensures
\[
\{ x : f(x) = \alpha \} ~\subset~ \Omega_{\e} ~\subset~ \{ x : f(x) \le \alpha \},
\]
so in fact  $\conv \Omega_{\e} = \{ x : f(x) \le \alpha \}$.  (Although superfluous for this proof, \mbox{\cite[Theorem 3.2]{yan}} even guarantees that $f_{\e}$ has positive-definite Hessian on this compact convex set, and hence is strongly convex on it.)

We next observe that the level set $\{ x \in U_{\e} : f_{\e}(x) \le f_{\e}(x_0) \}$ is bounded, since it is contained in the set $\{x : f(x) \le \alpha \}$.  Otherwise there would exist a point $x \in U_{\e}$ satisfying $f_{\e}(x) \le f_{\e}(x_0) = f(x_0) < \alpha$ and $f(x) > \alpha$.  By continuity of $f$, there exists a point $y$ on the line segment between $x_0$ and $x$ satisfying $f(y) = \alpha$.  But then we must have $y \in \Omega_{\e}$ and hence $f_{\e}(y) = f(y) = \alpha$, contradicting the convexity of $f_{\e}$.

The values and gradients of the functions $f$ and $f_{\e}\colon U_{\e} \to \R$ agree at each iterate 
$x_k$, so since those iterates comprise a BFGS sequence for $f$, they also do so for $f_{\e}$.  We can therefore apply Theorem \ref{thm1} to deduce 
\[
f(x_k) = f_{\e}(x_k) \downarrow \min f_{\e} ~~\mbox{as}~ k \to \infty.
\]
By assumption, there exists a sequence of points $x^r \in V$ (for $r=1,2,3,\ldots$) satisfying 
$\lim_r f(x^r) = \min f$.  For any fixed index $r$, we know $x^r \in \Omega_{\e}$ for all $\e > 0$ sufficiently small, so we have
\[
\min f \le \lim_k f(x_k) = \min f_{\e} \le f_{\e}(x^r) = f(x^r).
\]
Taking the limit as $r \to \infty$ shows $\lim_k f(x_k) = \min f$, as required.
\finpf

The following consequence suggests simple examples.

\begin{cor} \label{nonsmooth_min}
Powell's Theorem also holds with smoothness assumption (\ref{continuous}) replaced by the assumption that 
$\nabla^2 f$ is positive-definite and continuous throughout the set 
$\{ x \in U : f(x) > \min f \}$.
\end{cor}

\pf
Suppose the result fails.  The given set, which we denote $V$ must contain the set $\cl\!(x_k)$\/:  otherwise there would exist a subsequence of $(x_k)$ converging to a minimizer of $f$, and since the values $f(x_k)$ decrease monotonically, they would converge to $\min f$, a contradiction.  Clearly we have $\inf_V f = \min f$.  But now applying Theorem \ref{thm2} gives a contradiction.
\finpf

\begin{cor} \label{semi-algebraic}
Consider an open semi-algebraic convex set $U \subset \Rn$ containing a BFGS sequence for a semi-algebraic strongly convex function $f \colon U \to \R$ with bounded level sets.  Assume that the sequence and all its limit points lie in the interior of the set where $f$ is twice differentiable.  Then the sequence of function values converges to the minimum value of $f$.
\end{cor}

\pf
Denote the interior of the set where $f$ is twice differentiable by $V$.  Standard results in semi-algebraic geometry \cite[p.~502]{vandendries-miller} guarantee that $V$ is dense in $U$, whence $\inf_V f = \min f$, and furthermore that the Hessian $\nabla^2 f$ is continuous throughout $V$, and hence positive-definite by strong convexity.  The result now follows by Theorem~\ref{thm2}.
\finpf

The open set $V$ in the proof of Corollary \ref{semi-algebraic}, where the function $f$ is smooth, has full measure in the underlying set $U$.  Hence, if we initialize the algorithm in question with a starting point $x_0$ generated at random from a continuous probability distribution on $U$, and use a computationally realistic line search to generate each iterate $x_k$ from its predecessor, then we would expect $(x_k) \subset V$ almost surely.  Then, according to the result, one (or both) of two cases hold. 
\begin{enumerate}
\item[(i)]
The algorithm succeeds: $f(x_k) \to \min f$. 
\item[(ii)]
A subsequence of the iterates converges to a point where $f$ is nonsmooth.
\end{enumerate}
Extensive computational experiments suggest case (i) holds almost surely \cite{BFGS}.

Like Theorem \ref{thm2}, analogous versions of Corollary \ref{semi-algebraic} hold for many other algorithms, in addition to the BFGS method.  By contrast with BFGS, however, those algorithms often fail in general, due to the possibility of case (ii).  In the special situation described in Corollary \ref{nonsmooth_min}, case (ii) implies case (i), so analogous results will hold for many common algorithms, like steepest descent, coordinate descent, or conjugate gradients.

\section{Special constructions}
Unlike Powell's original result, Theorem \ref{thm2} requires the Hessian $\nabla^2 f$ to be positive-definite on an appropriate set, an assumption that fails for some simple but interesting examples like the  Euclidean norm.  We can sometimes circumvent this difficulty by a more direct construction, avoiding tools from \cite{yan}.  The following result is a version of Corollary \ref{nonsmooth_min} under a more complicated but weaker assumption.

\begin{thm} \label{thm3}
Powell's Theorem also holds with the smoothness assumption (\ref{continuous}) replaced by the following weaker condition:
\begin{quote}
For all constants $\delta > 0$, there is a convex open neighborhood $U_{\delta} \subset U$ of the set
$\{x \in U : f(x) \le f(x_0)\}$, and a  $C^{(2)}$ convex function 
$f_{\delta} \colon U_{\delta} \to \R$ satisfying $f_\delta(x) = f(x)$ whenever
$f(x_0) \ge f(x) \ge \min f + \delta$.
\end{quote}
\end{thm}

\pf
Clearly condition (\ref{continuous}) implies the given condition, since we could choose $U_{\delta} = U$ and $f_{\delta} = f$.  Assuming this new condition instead, suppose the conclusion of Powell's Theorem \ref{thm1} fails, so there exists a number $\delta > 0$ such that \mbox{$f(x_k) > \min f + 2\delta$} for all $k=0,1,2,\ldots$.  Consider the function $f_{\delta}$ guaranteed by our assumption.  Since $f$ is continuous, there exists a point $\bx \in U$ satisfying $f(\bx) = \min f + \delta$, and since $f_{\delta}(\bx) = f(\bx)$, we deduce 
$\min f_{\delta} \le \min f + \delta$.

Since $(x_k)$ is a BFGS sequence for the function $f$, it is also a BFGS sequence for the function $f_{\delta}$.  Applying Theorem \ref{thm1} with $f$ replaced by $f_{\delta}$ shows the contradiction
\[
\min f + 2\delta \le f(x_k) = f_{\delta}(x_k) \downarrow \min f_{\delta} \le \min f + \delta,
\]
so the result follows.
\finpf

We can apply this result directly to the Euclidean norm.

\begin{cor}
Any BFGS sequence for the Euclidean norm on $\Rn$ converges to zero.
\end{cor}

\pf
For any $\delta > 0$, consider the function $g_{\delta} \colon \R \to \R$ defined by
\bmye \label{g-delta}
g_{\delta}(t) = 
\left\{
\begin{array}{cl}
\frac{\delta^3 + 3\delta t^2 - |t|^3}{3\delta^2}  & (|t| \le \delta) \\ \\
|t| & (|t| \ge \delta).
\end{array}
\right.
\emye
This function is $C^{(2)}$ convex and symmetric.  The function $f_{\delta} \colon \Rn \to \R$ defined by $f_{\delta}(x)  = g_{\delta}(\|x\|)$ is also $C^{(2)}$ convex, either as a consequence of \cite{sendov_lorentz} or via a straightforward direct calculation.  The result now follows from Theorem \ref{thm3}.
\finpf

Analogously, the following result is a more direct version of Theorem \ref{thm2}.

\begin{thm} \label{thm4}
Powell's Theorem also holds with the smoothness assumption (\ref{continuous}) replaced by the assumption that some open set $V \subset U$ containing the set $\cl\!(x_k)$ and satisfying $\inf_V f = \min f$ also satisfies the following condition:
\begin{quote}
For all constants $\delta > 0$, there is a convex open neighborhood $U_{\delta} \subset U$ of the set
$\{x \in U : f(x) \le f(x_0)\}$, and a  $C^{(2)}$ convex function 
$f_{\delta} \colon U_{\delta} \to \R$ satisfying $f_\delta(x) = f(x)$ for all points $x \in U_{\delta}$ such that $d_{V^c}(x) > \delta$.
\end{quote}
\end{thm}

\pf
Denote the distance between the compact set $\cl\!(x_k)$ and the closed set $V^c$ by $\bar\delta$, so we know $\bar\delta > 0$.  For any constant $\delta \in (0,\bar\delta)$, we have $d_{V^c}(x_k) > \delta$ for all indices $k=0,1,2,\ldots$, and hence $f_\delta(x_k) = f(x_k)$.

The values and gradients of the functions $f$ and $f_{\delta}$ agree at each iterate 
$x_k$, so since those iterates comprise a BFGS sequence for $f$, they also do so for $f_{\delta}$.  We can therefore apply Theorem \ref{thm1} to deduce 
\[
f(x_k) = f_{\delta}(x_k) \downarrow \min f_{\delta} ~~\mbox{as}~ k \to \infty.
\]
By assumption, there exists a sequence of points $x^r \in V$ (for $r=1,2,3,\ldots$) satisfying 
$\lim_r f(x^r) = \min f$.  For any fixed index $r$, we know $d_{V^c}(x^r) > \delta$ for all sufficiently small $\delta > 0$, so since $f_{\delta}(x^r) = f(x^r)$, we deduce $\min f_{\delta} \le f(x^r)$.  The inequality $\lim_k f(x_k) \le f(x^r)$ follows, and letting $r \to \infty$ proves $\lim_k f(x_k) = \min f$ as required.
\finpf

We end by proving a claim from the introduction.

\begin{cor}
Any BFGS sequence for the function $f \colon \R^2 \to \R$ given by $f(u,v) = u^2 + |v|$ has a subsequence converging to a point on the line $v=0$.
\end{cor}

\pf
Suppose the result fails, so some BFGS sequence $\big((u_k,v_k)\big)$ has its closure contained in the open set
\[
V ~=~ \{ (u,v) \in \R^2 : v \ne 0 \}.
\]
Clearly we have $\inf_V f = \min f$.  For any constant $\delta > 0$, define a function 
$f_{\delta} \colon \R^2 \to \R$ by $f(u,v) = u^2 + g_{\delta}(v)$, where the function $g_{\delta}$ is given by equation (\ref{g-delta}).  Then we have $f(u,v) = f_{\delta}(u,v)$ for any point $(u,v)$ satisfying $|v|>\delta$, or equivalently $d_{V^c}(u,v) > \delta$.  Hence the assumptions of Theorem \ref{thm4} hold (using the set $U_\delta = \R^2$), so we deduce $f(u_k,v_k) \to 0$, and hence $(u_k,v_k) \to (0,0)$.  This contradiction completes the proof.
\finpf

As we remarked in the introduction, numerical evidence strongly supports a conjecture that all BFGS sequences for the function $f(u,v) = u^2 + |v|$ converge to zero.  That conjecture remains open.

\def\cprime{$'$} \def\cprime{$'$}


\end{document}